\documentclass[10pt]{article}

\usepackage{amsmath,amssymb,diagrams,newcent}
\usepackage{a4}

\newcounter{Prob}
\newenvironment{paragraaf}
{\stepcounter{Prob}%
\noindent \textbf{(\thesection.\arabic{Prob}) : \ }}{
\par \vskip 4mm}

\newcommand{\np}{{\par \noindent}}

\newcommand{\C}{\mathbb{C}}

\newcommand{\N}{\mathbb{N}}

\newcommand{\Z}{\mathbb{Z}}

\newcommand{\Oscr}{\mathcal{O}}

\newcommand{\atn}{\text{\bf @}_{\text{\bf n}}}

\newcommand{\kat}[1]{{\text{\bf \usefont{OT1}{pag}{b}{n} #1}}}
\newcommand{\wis}[1]{{\text{\bf \usefont{OT1}{pag}{m}{n} #1}}}

\begin{document}

\title{$\kat{noncommutative geometry}~\atn$}
\author{Lieven Le Bruyn\thanks{Research director of the FWO (Belgium)} 
\\ Universitaire Instelling Antwerpen \\
B-2610 Antwerp (Belgium)  \\
{\tt lebruyn$@\mbox{wins.uia.ac.be}$}}
\date{April 29, 1999}

\maketitle

\np
These are notes of a talk given in the 'noncommutative gometry' seminar at the
Max-Planck Instutute in Bonn, April 29th 1999. I thank Yu. Manin and Alex Rosenberg for
the invitation and the public for stimulating suggestions, some of which i have added
in footnotes.

\section{$\kat{smooth affine geometry}~\atn$.}

\begin{paragraaf} Let $\kat{cat}$ be a category of associative $\C$-algebras with unit. An
algebra $A \in Ob(\kat{cat})$ is said to be {\it $\kat{cat}$-smooth} iff for every {\it test-object}
$(B,I)$ in $\kat{cat}$ (that is, $B \in Ob(\kat{cat})$, $I \triangleleft B$ a nilpotent ideal such that
$\tfrac{B}{I} \in Ob(\kat{cat})$) and every $A \rTo^{\phi} \tfrac{B}{I} \in Mor(\kat{cat})$, there
exists a {\it lifting morphism} $A \rTo^{\tilde{\phi}} B \in Mor(\kat{cat})$ making the diagram below
commutative
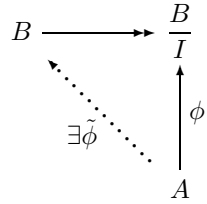
\[
\begin{diagram} B & \rOnto & \dfrac{B}{I} \\
& \luDotsto_{ \exists \tilde{\phi}} & \uTo_{\phi} \\
& & A \end{diagram}
\]
When $\kat{cat} = \kat{commalg}$, the category of all commutative $\C$-algebras, we recover Grothendieck's
formulation of smooth (regular) commutative algebras. For this reason we call a $\kat{commalg}$-smooth
algebra {\it g-smooth}.

\np
When $\kat{cat} = \kat{alg}$, the category of all $\C$-algebras, we recover Quillen's notion of
{\it quasi-free} or {\it formally smooth} algebras. For this reason we call a $\kat{alg}$-smooth
algebra {\it q-smooth}.

\np
Usually, we will assume that $\kat{cat}$-smooth algebras are {\it affine} algebras in $\kat{cat}$.

\np
Note however, that a commutative q-smooth algebra need not be g-smooth. For example, consider the
polynomial algebra $\C[x_1,\hdots,x_d]$ and the $4$-dimensional noncommutative local algebra
\[
B = \dfrac{\C \langle x,y \rangle}{(x^2,y^2,xy+yx)} = \C \oplus \C x \oplus \C y \oplus \C xy \]
Consider the one-dimensional nilpotent ideal $I = \C (xy-yx)$ of $B$, then the $3$-dimensional
quotient $\tfrac{B}{I}$ is commutative and we have a morphism $\C[x_1,\hdots,x_d] \rTo^{\phi} \tfrac{B}{I}$
by $x_1 \mapsto x, x_2 \mapsto y$ and $x_i \mapsto 0$ for $i \geq 2$. This morphism admits no lift
to $B$ as for any potential lift $[\tilde{\phi}(x),\tilde{\phi}(y)] \not= 0$ in $B$. Therefore,
$\C[x_1,\hdots,x_d]$ can only be q-smooth if $d=1$.
\end{paragraaf}

\begin{paragraaf} In fact, W. Schelter proved in \cite{Schelter} that $A$ is q-smooth if and only if
the $A$-bimodule $\Omega^1_A = Ker~A \otimes A \rTo^{m_A} A$ is projective. If $A$ is an
affine commutative $q$-smooth algebra, $A$ must be the coordinate ring of a finite set of points or
of a smooth affine curve. Apart from semisimple algebras and some algebras which are finite modules
over commutative $q$-smooth algebras, noncommutative q-smooth algebras are rather exotic objects
such as free algebras and algebras arising from universal constructions (some of which we will encounter
below).

\np
A fairly innocent class of q-smooth algebras are the {\it path algebras}. Let $Q$ be a {\it quiver},
that is a directed graph on a finite set $Q_v = \{ v_i,\hdots,v_k \}$ of vertices, having a finite
set $Q_a = \{ a_1, \hdots, a_l \}$ of arrows. The path algebra $\C~Q$ has as $\C$-basis the oriented
paths in $Q$ (including those of length zero corresponding to the vertices $v_i$) and multiplication
induced by concatenation, that is, $1 = v_1 + \hdots + v_k$ is a decomposition into
orthogonal idempotents and we have
\begin{itemize}
\item{$v_j.a$ is zero unless $\underset{v_j}{\bullet} \lTo^a \bullet$, }
\item{$a.v_i$ is zero unless $\bullet \lTo^a \underset{v_i}{\bullet}$, }
\item{$a_i.a_j$ is zero unless $\bullet \lTo^{a_i} \bullet \lTo^{a_j} \bullet$.}
\end{itemize}
To prove that $\C~Q$ is q-smooth, take a test-object $(B,I)$ in $\kat{alg}$ and an algebra map
$\C~Q \rTo^{\phi} \tfrac{B}{I}$. The decomposition $1 = \phi(v_1) + \hdots + \phi(v_k)$ into
orthogonal idempotents can be lifted modulo the nilpotent ideal to a decomposition
$1 = \tilde{\phi}(v_1) + \hdots + \tilde{\phi}(v_k)$ into orthogonal idempotents. But then, taking
for every arrow
\[
\underset{v_j}{\bullet} \lTo^a \underset{v_i}{\bullet} \quad \quad
\tilde{\phi}(a) \in \tilde{\phi}(v_j) (\phi(a) + I) \tilde{\phi}(v_i) \]
gives a required lift.
\end{paragraaf}

\begin{paragraaf} J. Cuntz and D. Quillen argue in \cite{CQ1,CQ2} that q-smooth algebras behave (for
example with respect to deRham cohomology) as commutative affine smooth algebras. In \cite[\S 9]{Kont1}
 M. Kontsevich gave a (somewhat cryptic) sketch how one might
go about to develop an affine noncommutative geometry for q-smooth algebras. He suggests
that one should {\it approximate} the noncommutative geometry of $A$ {\it at level $n$} by the
{\it representation space} $\wis{rep}_n~A$ which is the affine scheme representing the functor
\[
\kat{commalg} \rTo^{Hom_{\kat{alg}}(A,M_n(-))} \kat{sets} . \]
When $A$ is q-smooth, it follows from the lifting property for q- and g-smooth algebras that
$\wis{rep}_n~A$ is a smooth affine variety. We will denote this {\it approximation at level $n$} by
\[
\wis{space}~A \atn = \wis{rep}_n~A . \]
Recently, Kontsevich and A. Rosenberg \cite{KontRos}  made this proposal more
explicit. They argue that in order to extend a structure $\wis{struct}$ in commutative geometry to
the noncommutative affine smooth variety $\wis{spec}~A$ we must be able to define it at every level 
\[
\wis{struct}(\wis{space}~A) \Rightarrow \ \forall n \ : \ \wis{struct} \atn = \wis{struct}(\wis{rep}_n~A)
\]
That is, a noncommutative structure of some kind on $\wis{spec}~A$ should induce analogous
commutative structures on all the representation spaces $\wis{rep}_n~A$. These structures
$\wis{struct}$ can either be
\begin{itemize}
\item{{\it classical }, that is, ordinary commutative gadgets such as functions, vector fields and so
on (see \cite[\S 1.3.2]{KontRos} for more examples), or }
\item{{\it non-classical}, that is, new structures on commutative schemes coming from noncommutative
algebra such as the {\it formal structures} of Kapranov \cite{Kapranov} to be defined below.}
\end{itemize}
\end{paragraaf}

\begin{paragraaf} In this talk i want to indicate how one can extend this approximation proposal from
q-smooth algebras to arbitrary algebras and obtain in this way a rich
\wis{affine geometry} $\atn$. Then i will briefly indicate how one can construct global objects at
level $n$ and how one might build a noncommutative geometry from these approximate objects.
\end{paragraaf}

\section{$\kat{affine geometry}~\atn$.}

\setcounter{Prob}{0}

\begin{paragraaf} First we will specify $\kat{alg} \atn$, that is, the algebras that are level $n$
approximations of associative algebras.
A {\it trace map} on an associative $\C$-algebra $A$ is a $\C$-linear map $A \rTo^{tr} A$
such that for all $a,b \in A$ we have $tr(ab) = tr(ba)$, $tr(a)b = b tr(a)$ and $tr(tr(a)b) =
tr(a)tr(b)$. Algebras with trace are the objects of a category $\kat{alg}^{\kat{tr}}$
 with morphisms the trace preserving
$\C$-algebra maps. The forgetful functor $\kat{alg}^{\kat{tr}} \rTo \kat{alg}$ has a left adjoint
\[
\kat{alg} \rTo^{\tau} \kat{alg}^{\kat{tr}} \]
that is, given an algebra $A$ we can construct an algebra $A^{\tau}$ with trace in a universal way
by adding formally the traces.

\np
Fix a number $n$ and express $\prod_{i=1}^n (t-\lambda_i)$ as a polynomial in $t$ with coefficients
polynomials in the Newton functions $\nu_k = \sum_{i=1}^n \lambda_i^k$. Replacing $\nu_k$ by
$tr(x^k)$ we get a formal {\it Cayley-Hamilton polynomial of degree $n$} : $\chi^{(n)}_x(t)$. Let
$A$ be an algebra with trace $tr$, we say that $A$ is a {\it Cayley-Hamilton algebra of degree $n$}
if $tr(1) = n$ and $\chi^{(n)}_a(a) = 0$ in $A$ for all $a \in A$. Cayley-Hamilton algebras of
degree $n$ are the objects of a category $\kat{alg}^{\kat{tr}}_{\kat{n}}$ with trace preserving algebra
maps as morphisms. There is a natural functor
\[
\kat{alg}^{\kat{tr}} \rTo^n \kat{alg}^{\kat{tr}}_{\kat{n}} \]
by sending an algebra $A$ with trace $tr$ to the quotient $A_n$ by the twosided ideal generated by the
elements $tr(1)-n$ and $\chi_a^{(n)}(a)$ for all $a \in A$. Starting with an arbitrary algebra $A$
we propose to take $A~\atn = A^{\tau}_n$. That is,
\[
\begin{diagram}
\kat{alg} & & \\
\dTo^{\tau} \uInto & \rdDotsto^{\atn} & \\
\kat{alg}^{\kat{tr}} & \rTo_n & \kat{alg}^{\kat{tr}}_{\kat{n}} = \kat{alg} \atn
\end{diagram} \quad \quad \quad
\begin{diagram}
A & & \\
\dTo^{\tau} & \rdDotsto^{\atn} & \\
A^{\tau} & \rTo_n & A^{\tau}_n = A \atn
\end{diagram}
\]
Of course it may happen that $A~\atn = 0$ (for example if
$A = A_m(\C)$ the $m$-th Weyl algebra). For more details on algebras with trace and their properties
(some of which we will recall below) we refer to the paper \cite{ProcCH} of C. Procesi.
\end{paragraaf}

\begin{paragraaf} For $A$ an arbitrary associative $\C$-algebra we define $\wis{space}~A \atn$ as before
to be the representation space $\wis{rep}_n~A$ which is the affine scheme representing the functor
\[
\kat{commalg} \rTo^{Hom_{\kat{alg}}(A,M_n(-))} \kat{sets} \]
That is, there is a universal representation $A \rTo^{j_A} M_n(\C[\wis{rep}_n~A])$ such that for every $\C$-algebra
map $A \rTo^{\phi} M_n(B)$ with $B$ a commutative algebra, we have a uniquely determined morphism
$\C[\wis{rep}_n~A] \rTo^{\psi} B$ such that the diagram below is commutative
\[
\begin{diagram}
A & & \rTo^{\phi} & & M_n(B) \\
& \rdTo_{j_A} & & \ruTo_{M_n(\psi)} \\
& & M_n(\C[\wis{rep}_n~A]) & & 
\end{diagram}
\]
$GL_n$ acts by conjugation on $M_n(B)$ in a functorial way making $\wis{rep}_n~A$ into an affine $GL_n$-scheme
and so there are actions of $GL_n$ by automorphisms on $\C[\wis{rep}_n~A]$ and also on
$M_n(\C[\wis{rep}_n~A]) = \C[\wis{rep}_n~A] \otimes M_n(\C)$ (tensor product action). The image of
$A$ under $j_A$ is contained in the ring of invariants $M_n(\C[\wis{rep}_n~A])^{GL_n}$ which is the
ring of $GL_n$-equivariant maps $\wis{rep}_n~A \rTo M_n(\C)$ with algebra structure induced by the one
on the target space. 
The main
results on this $GL_n$-setting are due to Procesi \cite{ProcCH} and assert that we have a
commutative functorial diagram
\[
\begin{diagram}
A & \rTo^{j_A} & M_n(\C[\wis{rep}_n~A]) \\
\dTo^{\atn} & \rdTo^{j_A} & \uInto \\
A \atn & \cong & M_n(\C[\wis{rep}_n~A])^{GL_n}
\end{diagram}
\]
and we recover the algebra with trace $A \atn \in Ob(\kat{alg}^{\kat{tr}}_{\kat{n}})$ from the
$GL_n$-affine scheme $\wis{space}~A~\atn = \wis{rep}_n~A$. Clearly, this scheme is in general not
smooth nor even reduced. Procesi studied $\kat{alg}^{\kat{tr}}_{\kat{n}}$-smooth algebras in
\cite{ProcCH} which we will call p-smooth from now on. In fact, he proved that $A \atn$ is p-smooth
if and only if $\wis{rep}_n~A$ is a smooth $GL_n$-scheme. In particular, we have
\[
A \text{ q-smooth \ } \quad \Rightarrow \ \forall n \ : \quad A \atn \text{ p-smooth}  \]
consistent with our approximation strategy.
\end{paragraaf}

\begin{paragraaf}
We now come to classical structures on $\wis{space}~A~\atn = \wis{rep}_n~A$. Let us consider
functions as proposed in \cite[\S 1.3.2]{KontRos}. Let $a \in A$ and $A \rTo^x M_n(\C)$ a geometric
point $x \in \wis{rep}_n~A$, then we can take the ordinary trace $tr(x(a))$. This gives us a linear map
\[
\dfrac{A}{[A,A]} \rTo \C[\wis{rep}_n~A] \]
Kontsevich and Rosenberg propose to take for $\wis{funct}~A~\atn$ the subalgebra generated by the
image of this map. Observe that these functions are $GL_n$-invariant and recall that Procesi proved
that the ring of invariant polynomial functions
\[
\C[\wis{rep}_n~A]^{GL_n} = tr~A~\atn \]
where $tr$ is the trace map on $A~\atn$. That is, we have
\[
\wis{funct}~A~\atn = tr~A~\atn . \]
This ring has the following representation theoretic interpretation. The $GL_n$-orbits in
$\wis{rep}_n~A$ correspond to isomorphism classes of $n$-dimensional representations of $A$ and the
{\it closed} orbits correspond to $n$-dimensional {\it semi-simple} representations. Invariant theory
tells us that the closed orbits are parametrized by the maximal ideals of the ring of polynomial
invariants. The inclusion $tr~A~\atn \rInto \C[\wis{rep}_n~A]$ induces a morphism of schemes
\[
\wis{rep}_n~A \rOnto^{\pi} \wis{fac}_n~A \]
and sends an $n$-dimensional representation of $A$ to the isomorphism class of the
semi-simple $n$-dimensional representation which is the direct sum of the Jordan-H\"older factors.
Moreover, the algebra with trace $A~\atn$ is a finitely generated module over the subalgebra
$tr~A~\atn = \C[\wis{fac}_n~A]$ and hence we can associate to $A~\atn$ its classical structure
sheaf $\Oscr_{A~\atn}$ which is a sheaf of noncommutative algebras over $\wis{fac}_n~A$. 

\np
We propose that the classical structures on $\wis{space}~A~\atn$ are given by $GL_n$-equivariant
structures associated with the $GL_n$-invariant theoretic setting
\[
\begin{diagram}
& & \Oscr_{A \atn} \\
& & \dDots \\
\wis{rep}_n~A & \rOnto^{\pi} & \wis{fac}_n~A
\end{diagram}
\]
All the classical structures proposed in \cite[\S 1.3.2]{KontRos} are of this general form.
\end{paragraaf}

\begin{paragraaf}
Let us compute all of this in the case of path algebras. Consider the semisimple subalgebra
$V=\underbrace{\C \times \hdots \times \C}_k$ generated by the vertex-idempotents $\{ v_1,\hdots,v_k \}$.
Every $n$-dimensional representation of $V$ is semi-simple and determined by the multiplicities by
which the factors occur. That is, we have a decomposition
\[
\wis{rep}_n~V = \bigsqcup_{\sum a_i = n} GL_n / (GL_{a_1} \times \hdots \times GL_{a_k}) =
\bigsqcup_{\alpha} \wis{rep}_{\alpha}~V \]
into homogeneous spaces where $\alpha$ runs over the {\it dimension vectors} $\alpha = (a_1,\hdots,
a_k)$ such that $\sum_i a_i = n$. The inclusion $V \rInto \C~Q$ induces a map $\wis{rep}_n~\C~Q \rTo^{\psi}
\wis{rep}_n~V$ and we have the decomposition of $\wis{rep}_n~\C~Q$ into associated fiber bundles
\[
\psi^{-1}(\wis{rep}_{\alpha}~V) = GL_n \times^{GL(\alpha)} \wis{rep}_{\alpha}~Q . \]
Here, $GL(\alpha) = GL_{a_1} \times \hdots \times GL_{a_k}$ embedded along the diagonal in $GL_n$ and
$\wis{rep}_{\alpha}~Q$ is the affine space of $\alpha$-dimensional representations of the quiver
$Q$. That is,
\[
\wis{rep}_{\alpha}~Q = \bigoplus_{\underset{v_j}{\bullet} \lTo^a \underset{v_i}{\bullet}} 
M_{a_j \times a_i}(\C) \]
and $GL(\alpha)$ acts on this space via base-change in the vertex-spaces. That is, $\wis{rep}_n~\C~Q$
is the disjoint union of smooth affine components depending on the dimension vectors $\alpha = (a_1,
\hdots,a_k)$ such that $\sum a_i = n$. This
decomposition also translates to the algebra at level $n$
\[
\C~Q~\atn = \bigoplus_{\alpha}~T^Q_{\alpha} \]
where the trace map on the component $T^Q_{\alpha}$ is determined by $tr(v_i) = a_i$. For the
functions at level $n$ we have the decomposition
\[
tr~\C~Q~\atn = \bigoplus_{\alpha}~t^Q_{\alpha} \quad \text{where} \quad t^Q_{\alpha} = tr~T^Q_{\alpha}
\]
and is the ring of $GL(\alpha)$-invariants $\C[\wis{rep}_{\alpha}~Q]^{GL(\alpha)}$. In
\cite{LBProc} it was proved that this ring of invariants is generated by traces along oriented
cycles in the quiver $Q$ of length $\leq n^2$. That is, for fixed $\alpha$ replace any arrow
$\underset{v_j}{\bullet} \lTo^a \underset{v_i}{\bullet}$ by the $a_j \times a_i$ matrix of coordinate
functions on $\wis{rep}_{\alpha}~Q$ corresponding to the arrow $a$. Then, multiplying these matrices
along a cycle produces a square matrix whose trace is a $GL(\alpha)$-invariant.

\np
In \cite{LBProc} we also gave a computational description of the algebra with trace $T^Q_{\alpha}$. It
has a block decomposition
\[
T^Q_{\alpha} = \begin{bmatrix} T_{11} & \hdots & T_{1k} \\
\vdots & & \vdots \\
T_{k1} & \hdots & T_{kk} \end{bmatrix} \]
where $T_{ij}$ is the $t^Q_{\alpha}$-module spanned by the paths in the quiver $Q$ starting at vertex
$v_j$ and ending in vertex $v_i$ (again, the length of the required paths can be bounded by $n^2$).

\np
Observe that $t^Q_{\alpha}$ is positively graded and e will denote the $\mathfrak{m}$-adic
completion with respect to the graded maximal ideal $\mathfrak{m} = \oplus_{i \geq 1} t^Q_{\alpha}(i)$
by $\hat{t}^Q_{\alpha}$. Analogously we denote $\hat{T}^Q_{\alpha} = T^Q_{\alpha} \otimes \hat{t}^Q_{\alpha}$.
\end{paragraaf}

\begin{paragraaf}
In order for the approximation strategy for a q-smooth algebra $A$ to succeed we need to control the
local structure of $A~\atn$ and $\wis{fac}_n~A$. In fact a much more general result holds. Let $A$ be
an associative algebra and $x$ a geometric point of $\wis{fac}_n~A$. We say that $A~\atn$ is {\it
locally smooth} in $x$ provided $\wis{rep}_n~A$ is smooth along the closed orbit determined by $x$.
Observe that if $A$ is q-smooth (or even p-smooth), then $A~\atn$ is locally smooth in {\it all} $x$.

\np
The point $x$ determines an $n$-dimensional semi-simple representation $M_x$ of $A$, say
\[
M_x = S_1^{\oplus m_1} \oplus \hdots \oplus S_k^{\oplus m_k} \]
where $S_i$ is a simple $A$-representation of dimension $d_i$ occurring in $M$ with multiplicity
$m_i$, that is $n = \sum_i d_im_i$.

\np
We associate a {\it local quiver} $Q(x)$ to the point $x$. It has $k$ vertices $v_i$ (the number of
distinct simple components $S_i$ of $M_x$) and the number of arrows from $v_i$ to $v_j$ is given by
\[
\#~\underset{v_j}{\bullet} \lTo \underset{v_i}{\bullet} = \text{dim}_{\C}~Ext^1_A(S_i,S_j) . \]
Further, we define a {\it local dimension vector} $\alpha(x) = (m_1,\hdots,m_k)$ determined by the
multiplicities of the simple components in $M_x$.

\np
Denote the completion of the stalk of the structure sheaf $\Oscr_{\wis{fac}_n~A}$ in $x$ by
$\hat{\Oscr}_x$, then with notation as before we have an isomorphism of local algebras
\[
\hat{\Oscr}_x \cong \hat{t}^{Q(x)}_{\alpha(x)} . \]
This is an application of the {\it Luna slice theorem} in invariant theory, see \cite{LBsmooth}.
Further, if $\hat{\Oscr}^{A~\atn}_x$ denotes the algebra $A~\atn \otimes \hat{\Oscr}_x$, then it
is also proved in \cite{LBsmooth} that
\[
\hat{\Oscr}^{A~\atn}_x \underset{Morita}{\sim} \hat{T}^{Q(x)}_{\alpha(x)}  \]
where $\underset{Morita}{\sim}$ stands for {\it Morita-equivalence}, that is, equivalence of the
module categories. The precise form of the Morita equivalence is determined by the embedding
of $GL(\alpha(x)) = Stab_{GL_n}(M_x) \rInto GL_n$.

\np
Hence, if $A~\atn$ is locally smooth in $x$, the \'etale local structure of $\Oscr_{A~\atn}$ and
$\wis{fac}_n~A$ near $x$ is fully determined by the combinatorial data $Q(x)$ and $\alpha(x)$. A
lot more can be said about these combinatorics. For example, $(Q(x),\alpha(x))$ determines the
local quiver-data in neighbouring points, one can compute the dimensions of the strata in
$\wis{fac}_n~A$ consisting of points with the same local quiver-data and given the local dimension of
$\wis{fac}_n~A$ in $x$ one can even classify all the possible quiver-data. For some of these we 
refer to \cite{LBetale} and \cite{LBsmooth}.

\np
In particular, this applies to a q-smooth algebra $A$ in any point $x$ of $\wis{fac}_n~A$. That is,
the the approximation $A~\atn$ of {\it any} q-smooth algebra $A$ 
is locally (in the \'etale topology) isomorphic to the approximation of the subclass of path algebras
providing a handle on the exotic class of q-smooth affine algebras\footnote{A. Rosenberg suggests it
would be nice to
obtain the quiver-data from the non-commutative space as
defined in \cite[\S 2]{KontRos}. I will look into this.}
\end{paragraaf}

\begin{paragraaf}
Now, we turn to non-classical structures on $\wis{space}~A~\atn$ such as Kapranov's {\it formal
noncommutative structure}, see \cite{Kapranov}. 
Let $R$ be an associative
$\C$-algebra, $R^{Lie}$ its Lie structure and $R^{Lie}_m$ the subspace spanned by the
expressions $[r_1,[r_2,\hdots,[r_{m-1},r_m]\hdots ]$ containing $m-1$ instances of Lie brackets. The
{\it commutator filtration} of $R$ is the (increasing)
filtration by ideals $(F^d~R)_{d \in \Z}$ with $F^d~R = R$ for $d \in \N$
and
\[
F^{-d}~R = \underset{m}{\sum}~\underset{i_1 + \hdots + i_m = d}{\sum} R R^{Lie}_{i_1} R \hdots R R^{Lie}_{i_m} R 
\]
The associated graded $gr_F~R$ is a (negatively) graded commutative Poisson algebra with part of degree
zero $R_{ab} = \tfrac{R}{[R,R]}$. 

\np
Denote with $\kat{nil}_d$ the category of associative $\C$-algebras $R$ such that $F^{-d-1}R = 0$
(note that an algebra map is filtration preserving). Kapranov studied in \cite{Kapranov}
$\kat{nil}_d$-smooth algebras which we will call $k_d$-smooth algebras from now on.

\np
Kapranov proves \cite[Thm 1.6.1]{Kapranov} that any affine commutative smooth algebra $C$ has a
unique $k_d$-smooth {\it thickening} $R$ with $R_{ab} \simeq C$ (up to isomorphisms identical on $C$).
The approach of \cite[\S 1]{Kapranov} may be compared to that of M. Artin in \cite{ArtinDef}. 

\np
For
$R \in Ob(\kat{nil}_d)$, Kapranov introduces a sheaf $\Oscr_R$ of noncommutative
rings on the commutative scheme $X_{ab} = \wis{spec}~R_{ab}$. Observing that the commutator filtration for
such $R$ is {\it Zariskian} as in \cite{LiFVO} the approach of \cite[\S 2]{Kapranov} may be compared
to that of F. Van Oystaeyen in \cite[Chpt II]{FVObook}.

\np
If $X$ is an affine smooth commutative variety and $R_d$ the canonical $k_d$-smooth thickening of
$\C[X]$, then the sheaf of noncommutative algebras on $X$
\[
\Oscr^{f} = 
\underset{\leftarrow}{\text{lim}}~\Oscr_{R_d}
\]
is Kapranov's noncommutative formal structure on $X$. 

\np
If $A$ is q-smooth and affine, the representation space
$\wis{rep}_n~A$ is a smooth affine variety and hence is equipped with such a formal noncommutative
structure sheaf $\Oscr^f$. Part of the Kontsevich-Rosenberg proposal is that the level $n$
approximation $\wis{spec}~A~\atn$ should carry such a formal structure.
\end{paragraaf}

\begin{paragraaf} 
At first sight we face a serious problem as Kapranov's inductive construction of
$k_d$-smooth thickenings of a commutative $C$ only works when $C$ is g-smooth (in that case one
can define at each stage a {\it universal} central extension). Moreover, the construction of such a
formal structure must be sufficiently functorial to suit our purposes. 
I do not know of such a structure
 for all commutative affine $\C$-algebras (or even only those having a $GL_n$-action).

\np
Fortunately, we only need it for coordinate rings of representation spaces
$\C[\wis{rep}_n~A]$ and there we can apply some ringtheory of the early 70ties, in particular G. Bergman's
coproduct theorems \cite{Bergman} (see also the lecture notes of A. Schofield \cite[Chp. 2]{Schofield}
for more details). The starting point is that for every associative algebra $A$ the functor
\[
\kat{alg} \rTo^{Hom_{\kat{alg}}(A,M_n(-))} \kat{sets} \]
is {\it representable} in $\kat{alg}$. That is, there exists an associative $\C$-algebra
$\sqrt[n]{A}$ such that there is a natural equivalence between the functors
\[
Hom_{\kat{alg}}(A,M_n(-)) \underset{n.e.}{\sim} Hom_{\kat{alg}}(\sqrt[n]{A},-) . \]
In other words, for every associative $\C$-algebra $B$, there is a functorial one-to-one correspondence
between the sets
\[
\begin{cases}
\text{algebra maps} \quad A \rTo M_n(B)  \\
\text{algebra maps} \quad \sqrt[n]{A} \rTo B
\end{cases}
\]
To define $\sqrt[n]{A}$ consider the free algebra product $A \ast M_n(\C)$ and consider the subalgebra
\[
\sqrt[n]{A} = A \ast M_n(\C)^{M_n(\C)} = \{ p \in A \ast M_n(\C) \mid p.(1 \ast m) = 
(1 \ast m).p \ \forall m \in M_n(\C) \}
\]
Before we can prove the universal property of $\sqrt[n]{A}$ we need to recall a property that
$M_n(\C)$ shares with any Azumaya algebra : if $M_n(\C) \rTo^{\phi} R$ is an algebra morphism and
if $R^{M_n(\C)} = \{ r \in R \mid r.\phi(m) = \phi(m).r \ \forall m \in M_n(\C) \}$, then we have
$R \simeq M_n(\C) \otimes_{\C} R^{M_n(\C)}$. In particular, if we apply this to $R = A \ast M_n(\C)$
and the canonical map $M_n(\C) \rTo^{\phi} A \ast M_n(\C)$ where $\phi(m) = 1 \ast m$ we obtain that
$M_n(\sqrt[n]{A}) = M_n(\C) \otimes_{\C} \sqrt[n]{A} = A \ast M_n(\C)$.

\np
Hence, if $\sqrt[n]{A} \rTo^{f} B$ is an algebra map we can consider the composition 
\[
A \rTo^{id_A \ast 1} A \ast M_n(\C) \simeq M_n(\sqrt[n]{A}) \rTo^{M_n(f)} M_n(B) \]
to obtain an algebra map $A \rTo M_n(B)$. Conversely, consider an algebra map $A \rTo^g M_n(B)$ and
the canonical map $M_n(\C) \rTo^i M_n(B)$ which centralizes $B$ in $M_n(B)$. Then, by the universal
property of free algebra products we have an algebra map $A \ast M_n(\C) \rTo^{g \ast i} M_n(B)$ and
restricting to $\sqrt[n]{A}$ we see that this maps factors
\[
\begin{diagram}
A \ast M_n(\C) & \rTo^{g \ast i} & M_n(B) \\
\uInto & & \uInto \\
\sqrt[n]{A} & \rDotsto & B
\end{diagram} \]
and one verifies that these two operations are each others inverses. The algebra $\sqrt[n]{A}$ has other
surprising properties. For example, no matter how bad $A$ is, $\sqrt[n]{A}$ is a {\it domain} having
$\C^*$ as its group of invertible elements by \cite[Thm. 2.19]{Schofield}.
Now, equip $\sqrt[n]{A}$ with the commutator filtration
\[
\begin{diagram}
\hdots & & -2 & & -1 & & 0 & & 1 &  & \hdots \\
\hdots & \subset & F_{-2} & \subset & F_{-1} & \subset & \sqrt[n]{A} & =  & 
\sqrt[n]{A} & = & \hdots 
\end{diagram}
\]
By the universal property of $\sqrt[n]{A}$ we see that $\tfrac{A}{F_{-d-1}} \in Ob(\kat{nil}_d)$ is the
object representing the functor
\[
\kat{nil}_d \rTo^{Hom_{\kat{alg}}(A,M_n(-))} \kat{sets} . \]
In particular, as $\kat{nil}_0 = \kat{commalg}$ we deduce that
\[
gr_0~\sqrt[n]{A} = \dfrac{\sqrt[n]{A}}{F_{-1}} = \dfrac{\sqrt[n]{A}}{[\sqrt[n]{A},\sqrt[n]{A}]} \simeq
\C[\wis{rep}_n~A] \]
as both algebras represent the same functor. The construction of the formal sheaf of
noncommutative algebras is similar to that of \cite{Kapranov}. For fixed $d$, the induced
filtration on the quotient $\sqrt[n]{A}_d = \tfrac{\sqrt[n]{A}}{F_{-d-1}}$ is Zariskian (even discrete). Hence,
taking the {\it saturation} of a multiplicatively closed subset of the associated graded
$gr~\sqrt[n]{A}$ is an Ore-set in $\sqrt[n]{A}_d$. Hence, we can construct a sheaf of algebras
$\Oscr_{\sqrt[n]{A}_d}$ on $\wis{rep}_n~A$ as in \cite{FVObook}. The formal noncommutative structure sheaf on
$\wis{rep}_n~A$ is then the
inverse limit
\[
\Oscr_{\sqrt[n]{A}}^f = \underset{\leftarrow}{lim}~\Oscr_{\sqrt[n]{A}_d} \]
We now claim that the approximation at level $n$ of $\wis{spec}~A$ is given by the {\it vier-span}
\[
\wis{spec}~A~\atn \quad = \quad \quad \begin{diagram}
\Oscr_{\sqrt[n]{A}}^f & & \Oscr_{A \atn} \\
\dDots & & \dDots \\
\wis{rep}_n~A & \rOnto & \wis{fac}_n~A
\end{diagram}
\]
There remains to prove that the formal structure $\Oscr_{\sqrt[n]{A}}^f$ defined above coincides
with Kapranov's formal structure in case $A$ is q-smooth. We have seen that
\[
M_n(\sqrt[n]{A}) \simeq A \ast M_n(\C) \]
hence if $A$ is q-smooth, so is $\sqrt[n]{A}$ by the coproduct theorems, see for example
\cite[Thm. 2.20]{Schofield} for a strong version. But then clearly $\sqrt[n]{A}_d$ is $k_d$-smooth
for all $d$ (or use \cite[Prop. 1.4.6]{Kapranov}). By the unicity of $k_d$-smooth thickenings
\cite[\S 1.6]{Kapranov} it is then immediate that Kapranov's formal structure $\Oscr^f$ on
$\wis{rep}_n~A$ coincides with our $\Oscr_{\sqrt[n]{A}}^f$.

\np
Note also the perhaps surprising fact that for $\wis{rep}_n~A$ for $A$ a q-smooth algebra one does
not need the perturbative approach of Kapranov to describe a formal neighborhood of
$\wis{rep}_n~A$ into a noncommutative smooth space. By the above, $\wis{rep}_n~A$ is a closed
subvariety of the noncommutative smooth space $\wis{spec}~\sqrt[n]{A}$.

\np
Our description also clarifies the importance of this formal noncommutative structure. In order
to understand the (commutative) scheme structure of $\wis{rep}_n~A$ one needs to consider also
representation $A \rTo M_n(F)$ where $F$ is a finite dimensional commutative algebra (rather than
restrict to $\C$). Similarly, if one wants to understand the noncommutative scheme structure of
$\wis{rep}_n~A$ one ought to consider representations $A \rTo M_n(F)$ where $F$ is a finite
dimensional non-commutative algebra. If $F$ is {\it basic} (that is, all simple $F$-representations
are one-dimensional over $\C$), then $F \in Ob(\kat{nil}_d)$ for some $d$, and then this representation
is controlled by the formal structure. In general, any finite dimensional noncommutative algebra
$F$ is Morita equivalent to a basic algebra, so certainly if we vary $n$ all the formal structures
$\Oscr_{\sqrt[n]{A}}$ control representation $A \rTo M_m(F)$. Expressed differently, one can view
the formal structures $\Oscr_{\sqrt[n]{A}}$ as sewing-machines to stitch the different
$\wis{rep}_n~A$ together.
\end{paragraaf}

\begin{paragraaf} I cannot resist the temptation to add an infinite family of formal structures on
$\wis{space}~A~\atn$. If at layer $1$ the commutator filtration is defined using the Lie bracket
$[r_1,r_2]$, one can similarly define at layer $m$ the $m$-commutator filtration based on the
expressions
\[
S_{2m}(r_1,r_2,\hdots,r_{2m}) = \sum_{\sigma \in S_{2m}} (-1)^{sgn(\sigma)} r_{\sigma(1)}
r_{\sigma(2)} \hdots r_{\sigma(2m)} . \]
The associated graded algebra with respect to the $m$-commutator filtration on $\sqrt[n]{A}$
is no longer
commutative but is a {\it polynomial identity algebra} of degree $m$ (that is, it basically lives
at level $\text{\bf @}_{\text{\bf m}}$). As these algebras are close to commutative algebras, they
have plenty of Ore sets and again the Zariskian argument provides us with new formal structure
sheaves.
\end{paragraaf}

\section{$\kat{geometry}~\atn$.}

\setcounter{Prob}{0}

\begin{paragraaf} In this section we will briefly indicate how one can define global objects
in the approximate geometry $\atn$ and afterwards how one can put approximate objects together
to define a noncommutative geometry. More details will have to await another occasion.

\np
The strategy to define a $\kat{geometry}~\atn$ is simple : use $GL_n$-equivariant (commutative)
geometry as inspiration and try to define all concepts in terms of $\kat{alg}~\atn$. The latter is
essential in order to define formal structures as we have seen above.

\np
In view of our local combinatorial description of p-smooth algebras, it is clear that the natural
topology of $\kat{affine geometry}~\atn$ is the \'etale topology. Therefore, in order to define
global objects we choose for the approach via {\it algebraic spaces}\footnote{V. Hinich suggests
it might be more natural to extend the notion of algebraic stack in order to maintain
compatibility with the Kontsevich-Rosenberg proposal of noncommutative spaces in
\cite[\S 2]{KontRos}. I agree and will try to work this out.} as developed by M. Artin
\cite{ArtinAS}. Algebraic spaces are defined by \'etale equivalence relations, so we need
\begin{itemize}
\item{a product in $\kat{alg}~\atn$, and}
\item{\'etale morphisms in $\kat{alg}~\atn$.}
\end{itemize}
As a product, the tensor-product $\otimes_{\C}$ is not suitable as it takes us out of
$\kat{alg}~\atn$ (think of tensorproducts of matrixrings). However, free algebra products
provide us with a suitable definition
\[
A \underset{n}{\boxtimes} B \overset{def}{=} A \ast B~\atn
\]
In order to define \'etale morphisms in $\kat{alg}~\atn$ we extend the notions of formally
\'etale (resp. formally unramified, formally smooth) from $\kat{commalg}$ verbatim. That is,
let $A$ be an associative algebra, then a morphism $A~\atn \rTo^{\phi} B$ in $\kat{alg}~\atn$
is said to be {\it formally \'etale} iff for every test-object $(T,I)$ in $\kat{alg}~\atn$,
we have a unique lift, where all the morphisms below are $A~\atn$-algebra maps
\[
\begin{diagram}
T & \rOnto & \dfrac{T}{I} \\
\uTo & \luDotsto_{\exists ! \tilde{\psi}} & \uTo^{\psi} \\
A~\atn & \rTo_{\phi} & B 
\end{diagram}
\]
If one replaces unicity by existence $\phi$ is said to be {\it formally smooth} and if unicity is
replaced by the existence of at most one lift, $\phi$ is said to be {\it formally unramified}. These
notions also have a geometric interpretation :  consider an algebra map
 $A \rTo^{\phi} B$, then $A~\atn \rTo^{\phi~\atn} B~\atn$ is \'etale if and only if the induced
 morphism $\wis{rep}_n~B \rTo \wis{rep}_n~A$ is an \'etale morphism of commutative schemes.

\np
In order to show that this {\it \'etale topology} on $\wis{spec}~A~\atn$ is rich enough, let us give
tow classes of \'etale maps. The first is classical, coming from commutative \'etale maps. Consider the
diagram
\[
\begin{diagram}
A~\atn & \rTo^{id \otimes f} & A~\atn \otimes_{tr~A~\atn} S \\
\uInto & & \uTo \\
tr~A~\atn & \rTo^f & S
\end{diagram}
\]
where $f$ is an \'etale morphism in $\kat{commalg}$, then $id \otimes f$ is \'etale in $\kat{alg}~\atn$.
The second class is more exotic, though it is the natural substitute for localizations, {\it
universal localizations} in $\kat{alg}~\atn$. Let $\kat{projmod}~A$ denote the category of finitely
generated (left) modules over an associative algebra $A$ and let $\Sigma$ be some class of maps in
this category (that is some $A$-module morphisms between certain projective modules). In \cite[Chp. 4]{Schofield}
it is shown that there exists an algebra map $A \rTo^{j_{\Sigma}} A_{\Sigma}$ with the universal
property that the maps $A_{\Sigma} \otimes_A \sigma$ have an inverse for all $\sigma \in \Sigma$.
Using the above geometric characterization of \'etale maps it follows that the induced maps
\[
A~\atn \rTo^{j_{\Sigma}~\atn} A_{\Sigma}~\atn \]
are \'etale in $\kat{alg}~\atn$. In particular, it may happen that a finite dimensional algebra has
hugely infinite dimensional \'etale extensions in $\kat{alg}~\atn$. For example consider the
$m+3$-dimensional algebra
\[
A = \begin{bmatrix} \C & V \\ 0 & \C \end{bmatrix} \]
where $V$ is a $\C$-vectorspace of dimension $m+1$. Then, there is a universal localization of $A$
isomorphic to $M_2(\C \langle x_1,\hdots,x_m \rangle)$. We refer to the book \cite{Schofield} of
A. Schofield for more details.

\np
Using these \'etale maps and the product $\underset{n}{\boxtimes}$ one can define
$\kat{algebraic spaces}~\atn$ as in \cite{ArtinAS}.
\end{paragraaf}

\begin{paragraaf} It may be interesting to reconsider the notion of {\it orbit topos} introduced by
R.W. Thomason in \cite{Thomason} in this setting. For example, one wonders how much of the orbit topos
of $\wis{rep}_n~A$ can be described using algebraic spaces in $\kat{alg}~\atn$.

\np
In general it is not true that any $GL_n$-equivariant technique can be performed in 
$\kat{algebraic spaces}~\atn$.
A noteworthy exception is $GL_n$-equivariant desingularization. In \cite{LBetale} the
following problem was handled. Let $X$ be a smooth (commutative) surface and $\Delta$ a central
simple algebra in $\kat{alg}~\atn$ over the function field $\C(X)$. Recall that $\Delta$ is determined by
\begin{itemize}
\item{a divisor $D \rInto X$ and a list of its irreducible components $C$,}
\item{the list of singular points $p_j$ of $D$, and}
\item{for each branch $B_k$ of $D$ at $p_i$ a number $n_{i,k} \in \Z/n\Z$ such that
$\sum_k n_{i,k} = 0$.}
\end{itemize}
One might ask whether for each $\Delta$ there exists a smooth object in $\kat{algebraic spaces}~\atn$
having $\Delta$ as its noncommutative function algebra. An idea might be to start with a {\it
maximal order} $\mathcal{A}$ in $D$ over $X$, consider $\wis{rep}_n~\mathcal{A}$ and construct its
$GL_n$-equivariant desingularization. In \cite[Chp. 6]{LBetale} it was shown that the above problem
has a positive solution if and only if
\[
\sum_{p_i \in C} \sum_{B_k \in C} n_{i,k} = 0 \quad \text{ for all irreducible components $C$ of $D$}
\]
In general, one can construct an object in $\kat{algebraic spaces}~\atn$ having $D$ as function algebra
and having only a finite number of points where it is not locally smooth. All of these singularities
are locally (in the \'etale topology) Morita equivalent to that of the quantum plane
$\C_q[x,y]$ with $xy = q yx$ for $q$ an $m$-th root of unity where $m \mid n$.

\np
Based on this example, our strategy to develop $\kat{geometry}~\atn$ is to use $GL_n$-equivariant
theory as far as it can serve us, but then study the remaining cases (which will lead to interesting
algebras) rather than solving the remaining problems by extending our algebraic framework.
\end{paragraaf}

\begin{paragraaf}
Finally, how can one use $\kat{geometry}~\atn$ to develop a noncommutative geometry. Again the strategy
is simple : formalize the setting $A \Rightarrow~\forall n~A~\atn$ and make a list of relations
holding naturally among the $\wis{spec}~A~\atn$. Then, define an object in $\kat{geometry}$ by a
list of objects $X~\atn$ in $\kat{geometry}~\atn$ satisfying this list of relations. At this
moment i do not have an elegant set of axioms. Let me conclude by proposing one axiom which
should illustrate the principle.

\np
If $n = \sum_i m_i$ and let $V_i$ be an $m_i$-dimensional representation of $A$, then the direct
sum $\oplus V_i$ is an $n$-dimensional representation of $A$. Hence there are morphisms
\[
\underset{i}{\times}~\wis{spec}~A~\text{\bf @}_{{\text{\bf m}}_i} \rTo \wis{spec}~A~\atn \]
satisfying obvious compatibility relations.

\np
To formalize this condition for objects $X = (X~\atn)_n$ in $\kat{geometry}$, let us stratify the
geometric points of the underlying $GL_n$-algebraic space $X~\atn$ by
\[
X~\atn(r) = \{ x \in X~\atn \mid Stab_{GL_n}(x) \ \text{has a maximal torus of dimension $r$} \ \} \]
Then, for each integral solution $n = \sum_i m_i$ one must have connecting morphisms
\[
\underset{i}{\times} X~\text{\bf @}_{{\text{\bf m}}_i} \rTo^{c_{(m_i)}} X~\atn \]
with the additional condition that
\[
X~\atn(r) \subset \underset{m_1+\hdots+m_r=n}{\cup} GL_n.Im(c_{(m_1,\hdots,m_r)}). \]
It is clear that a lot of additional work needs to be done.
\end{paragraaf}

\end{document}